\documentclass[11pt,bezier]{article}
\usepackage{amsmath,amssymb,amsfonts,euscript}
\usepackage{algorithmic, algorithm}
\usepackage{graphicx}

\usepackage[ansinew]{inputenc}
\usepackage{graphicx}
\usepackage{color}
\usepackage{mathrsfs}

\usepackage[colorlinks]{hyperref}
\usepackage[active,new,noold,marker]{xrcs}
\catcode`\=13
\def{$\bowtie$}

\usepackage{eurosym}
\usepackage{lscape} 

\newtheorem{prethm}{{\bf Theorem}}

\newenvironment{thm}{\begin{prethm}{\hspace{-0.5
               em}{\bf}}}{\end{prethm}}

\newtheorem{prepro}{{\bf Theorem}}

\newtheorem{preprop}{{\bf Proposition}}

\newtheorem{precor}{{\bf Corollary}}

\newenvironment{cor}{\begin{precor}{\hspace{-0.5
               em}{\bf}}}{\end{precor}}

\newtheorem{preconj}{{\bf Conjecture}}

\newtheorem{predefi}{{\bf Definition}}

\newtheorem{preremark}{{\bf Remark}}

\newenvironment{remark}{\begin{preremark}\rm{\hspace{-0.5
               em}{\bf}}}{\end{preremark}}

\newtheorem{preexample}{{\bf Example}}

\newenvironment{example}{\begin{preexample}\rm{\hspace{-0.5
               em}{\bf}}}{\end{preexample}}

\newtheorem{prelem}{{\bf Lemma}}

\newtheorem{prelam}{{\bf Lemma}}

\newtheorem{preprob}{{\bf Problem}}

\newtheorem{preproof}{{\bf Proof}}

\newtheorem{preali}{{\bf Proof of Theorem 1.}}

\newenvironment{ali}[1]{\begin{preali}{\rm
               #1}\hfill{$\Box$}}{\end{preali}}

\newtheorem{prealii}{{\bf Proof of Theorem 2.}}

\newenvironment{alii}[1]{\begin{prealii}{\rm
               #1}\hfill{$\Box$}}{\end{prealii}}

\newtheorem{prealiii}{{\bf Proof of Theorem 3.}}

\newenvironment{aliii}[1]{\begin{prealiii}{\rm
               #1}\hfill{$\Box$}}{\end{prealiii}}

\title{ On the complexity of deciding whether the regular
number is at most two}

\author{{\normalsize
{\sc Ali Dehghan${}^{\mathsf{a}}$},\,
{\sc Mohammad-Reza Sadeghi${}^{\mathsf{a}}$},\,
 {\sc Arash Ahadi${}^{\mathsf{b}}$}\,
}\vspace{3mm}
\\{\footnotesize{${}^{\mathsf{a}}$\it Department of Mathematics and Computer Science,
Amirkabir University of Technology, Tehran,
Iran}}  {\footnotesize{}}\\{\footnotesize{${}^{\mathsf{b}}$\it
Department of
Mathematical Sciences, Sharif University of Technology, Tehran, Iran}}
\thanks{{\it E-mail addresses}:  $\mathsf{ali\_dehghan16@aut.ac.ir}$, $\mathsf{msadeghi@aut.ac.ir}$, $\mathsf{arash\_ahadi@mehr.sharif.edu}$. } }

\date{}
\begin{document}
\maketitle

\begin{abstract}
{\small \noindent

The {\it regular number} of a graph $G$ denoted by $reg(G)$ is the minimum number of subsets into which the edge set of $G$ can be partitioned
so that the subgraph induced by each subset is regular.
In this
work we answer to
the problem posed
as an
open problem in A. Ganesan et al. (2012) \cite{kulli2}
about the complexity of determining the regular
number of  graphs.
We show that computation of the regular number for connected bipartite graphs is {\bf NP}-hard.
Furthermore, we show that, determining whether $ reg(G) = 2 $ for a given connected $3$-colorable graph $ G $  is {\bf NP}-complete. Also, we prove that a
new variant of the Monotone Not-All-Equal 3-Sat problem is {\bf NP}-complete.
}

\begin{flushleft}


\end{flushleft}

\end{abstract}

\section{Introduction}
\label{}
All graphs considered here are finite, undirected and loopless.
For every $v\in V (G)$, $d(v)$ denotes the degree of $v$ and for a natural number $k$, a graph $G$ is
called a {\it $k$-regular graph } if $d(v) = k$, for each $v \in V(G)$. Also, we denote the maximum degree
 of $G$ by $\Delta(G)$. We
follow \cite{MR1367739} for terminology and notation where they are not defined here.

In 1981,
Holyer \cite{MR635429}
considered the computational complexity of edge-partition problems and
 proved that for each fixed $n  \geq 3$ it is {\bf NP}-complete to determine whether an arbitrary graph
can be edge-partitioned into subgraphs isomorphic to the complete graph $K_n$.
Holyer purposed the following problem and conjectured that for every graph $H$ with more than two edges the following problem is {\bf NP}-complete.

{\em The Holyer Problem}:
Given a fixed graph $H$, can the edge set of an instance graph $G$ be partitioned into
subsets inducing graphs isomorphic to $H$?

Afterwards, the complexity of edge-partition problems have been studied extensively by several
authors, for instance see \cite{MR2510539,MR1460720,MR635429,  MR2297163}.
Dor et al. proved that the Holyer problem
is {\bf NP}-complete whenever $H$ has a connected component with at least 3 edges \cite{MR1460720}.
Afterwards, Bry{\'s} et al. in \cite{MR2510539} completely solved the Holyer problem;
they proved that, if a graph $H$ does not have a component with at least $3$ edges then the problem of deciding if the edge set of an instance graph $G$ can be partitioned into subsets inducing graphs isomorphic to
$H$  is polynomial time solvable.
Nowadays, computational complexity of edge decomposition problems is a well-studied area of graph theory.

If we consider the Holyer problem for a family $\mathcal{Q}$ of graphs instead of $H$ as a fixed graph
then, we can find some interesting problems.
For a family $\mathcal{Q}$ of graphs, a  $\mathcal{Q}$-decomposition of a graph $G$ is a partition of the edge set of $G$ into
subgraphs isomorphic to members of $\mathcal{Q}$.
Problems of $\mathcal{Q}$-decomposition of graphs have received a considerable attention,
for example, Holyer proved that it is {\bf NP}-hard to edge-partition a graph into the minimum number of complete subgraphs \cite{MR635429}. For more examples see \cite{MR1403227, MR679633}.

Consider the family of regular graphs, the edge set of every graph can be partitioned
such that the subgraph induced by each subset is regular.
In 2001, Kulli et al. introduced an interesting parameter for the edge-partition of a graph \cite{kulli}.
The {\it regular number} of a graph $G$ denoted by $reg(G)$ is the minimum number of subsets into which the edge set of $G$ can be partitioned
so that the subgraph induced by each subset is regular.
Nonempty subsets $E_1, \ldots ,E_r$ of $E$ are said to form a regular
partition of $G$ if the subgraph induced by each subset is regular.

\begin{example}{
The zebra $B_n$ is a graph defined recursively.
The zebra $B_0$ consists of a copy of $K_1$. In order to define $B_{n+1}$,
 consider a copy of $B_n$ and call it $I_n$, Also consider $|V(B_n) |=3^n$ new isolated vertices and call them $J_n$. Now put $2^n$
distinct prefect matchings between $I_n$  and $J_n$ call the resulting graph $B'_n$. Let $B_{n+1}=B'_n \cup B_n$.  The set of degrees of $B_n$ is $ \{0, 1,2,\cdots,2^n-1\}$
and it is easy to see that $reg(B_n)= n$.
}\end{example}

The  edge chromatic number
of a
graph denoted by $\chi '(G)$ is the minimum size of a partition of the edge set into $1$-regular subgraphs.
By Vizing's theorem \cite{MR0180505}, the edge chromatic number of a graph $G$ is equal to either $ \Delta(G) $ or $ \Delta(G) +1 $.
Therefore the regular number problem is a generalization for the  edge chromatic number and we have the following.

\begin{equation}
 reg(G)\leq \chi '(G) \leq \Delta(G) +1.
\end{equation}

In this
work we give the answer to
the problem posed
as an
open problem in A. Ganesan et al. \cite{kulli2}
about the complexity of determining the regular
number of a graph.
Determining the computational complexity of the regular number for disconnected graphs is a very easy problem. In the following we first consider disconnected graphs and then we focus on  the complexity boundary of regular number problem for connected graphs.

\begin{remark}{ It was shown that it is $ \mathbf{NP} $-hard to determine the edge chromatic number of an
$r$-regular graph; for any $r \geq 3$ \cite{MR689264}.
For a given $r$-regular graph $G$ with $r \geq 3$, consider the new graph $G'=G \cup K_{1,r}$. Clearly, $reg(G')=\Delta(G') $ if and only if $\chi'(G)=\Delta(G)$.
Therefore,
for every $k\geq 3$, it is
{\bf NP}-complete to decide whether $reg(G)=k$ for a given disconnected graph
$G$.
}\end{remark}

\begin{cor}
For every $k\geq 3$, it is
{\bf NP}-complete to decide whether $reg(G)=k$ for a given disconnected graph
$G$.
\end{cor}

In every regular partition of a tree $T$, each part is $1$-regular, so $reg(T)=\Delta(T)=\chi '(G)$ \cite{kulli2}.
Designing an algorithm to decompose a given bipartite graph into the minimum number
of regular subgraphs was posed as a problem in \cite{kulli2}.
We show that computation of the regular number is {\bf NP}-hard for connected bipartite graphs.

\begin{thm}
Computation of the regular number for connected bipartite graphs is {\bf NP}-hard.
\end{thm}

For a connected graph $G$, $reg(G)=1$, if and only if $G$ is regular. We consider the complexity of deciding whether $ reg(G) = 2 $ for a given graph $G$.
We prove that if $  \{ d(v): v\in V(G)\}  =\{k,k'\}$, where $k$ and $ k'$ are different natural numbers, then computation of the regular number is {\bf NP}-hard even if the graph is $3$-colorable.

\begin{thm}\label{T1}
Determining whether $ reg(G) = 2 $ for a given connected $3$-colorable graph $ G $  is {\bf NP}-complete.
\end{thm}

It
was asked  in \cite{kulli2} to determine whether  $reg(G)\leq \Delta(G)$ holds for all
connected graphs $G$.
We show that not only  there exists a counterexample for the above bound but also for a given connected graph $G$ decide whether $reg(G)=\Delta(G)+1$ is
{\bf NP}-complete.

\begin{thm}\label{T3}
Determining whether $reg(G)\leq \Delta(G)$ for a given connected graph $ G $  is {\bf NP}-complete.
\end{thm}

\section{NP-Completeness}

First, we show that computation of the regular number is {\bf NP}-hard for connected bipartite graphs.
It is shown that {\em 3-Partition} is $ \mathbf{NP}$-complete in the strong sense \cite{MR1567289}.

{\em 3-Partition.}\\
\textsc{Instance}: A positive integer $ k \in \mathbb{Z}^{+}$ and $3n$ positive integers $a_1, \ldots, a_{3n} \in \mathbb{Z}^{+}$ such that $  k/4 < a_i < k/2$ for each $1 \leq i \leq 3n$ and $ \sum _{i=1}^{3n} a_i=nk$.\\
\textsc{Question}: Can $\{a_1, \ldots, a_{3n}\}$ be partitioned into $n$ disjoint sets $A_1, \cdots , A_n$ such that,
for $1 \leq i \leq n$, $ \sum _{a \in A_i}a= k$?\\

\begin{ali}{
We reduce {\em 3-Partition} to
our problem in polynomial time.
For an instance $A=<a_1, \ldots,a_{3n}>$ and number $k$, consider $3n$ copies of complete bipartite graph $K_{2k,2k-1}$ and denote them by $K^{(1)}, \cdots, K^{(3n)}$. Also, denote the part of size $2k$ of $K^{(i)}$, by $X^{(i)}$. Next, consider three new vertices $u$, $v$ and $w$. For every $i$, $1 \leq i \leq 3n$, join $2a_i$ vertices of $X^{(i)}$ to $u$. From other vertices of $X^{(i)}$, join $k-a_i$ vertices of $X^{(i)}$ to $v$ and finally join $k-a_i$ remaining vertices of $X^{(i)}$ to $w$.
Call the set of these $2k$ edges $E^{(i)}$.
Also, call the resulting graph $G=G_{A,k}$.
We show that $reg(G) \leq n $ if and only if {\em 3-Partition} has a true answer. Since the degrees of $u$, $v$ and $w$ are $2kn$ and the degrees of all other vertices is $2k$, so in every regular partition $G_1, \cdots, G_t$ of $G$, each $G_i$ is at most $(2k)$-regular, therefore, $reg(G) \geq n $.

First, assume that $reg(G)=n$ and $G_1, \cdots, G_n$ is a regular partition of $G$ such that $G_i$ is $(r_i)$-regular. Since in every regular partition $G_1, \cdots, G_t$ of $G$, each $G_i$ is at most $(2k)$-regular, so $r_1=\cdots = r_n = 2k$. Hence, for every $i$,
$1 \leq i \leq  3n$, all the edges of $E^{(i)}$ and $K^{(i)}$ are in one of the $G_1, \cdots, G_n$. Consequently, for every $j$, we have $\sum_{i : K^{(i)}\subseteq G_j}2a_i=2k$. So $\{ \{ a_i : K^{(i)}\subseteq G_j\} : 1 \leq j \leq n \}$ is an appropriate partition for $A=<a_1, \ldots,a_{3n}>$.
On the other hand, if $A_1, \cdots, A_n$ is a partition for $A$ such that  $ \sum _{a \in A_i}a= k$, then for every $j$, let $G_j$ be the induced subgraph on $\{ u,v,w\}\bigcup_{i:a_i\in A_j} K^{(i)}$. So  $reg(G) \leq n $.
}\end{ali}

\begin{alii}{

Our proof consist of two steps. In the first step we prove that a new variant of the Monotone Not-All-Equal 3-Sat
problem is $ \mathbf{NP}$-complete, then in the step 2 we reduce it to our problem in polynomial time.

{{\bf Step 1.}} Let $\Upsilon$ be a $3$-SAT formula with the set
of clauses $C$ and the set
of variables
$X $. It is shown that the following problem is $ \mathbf{NP}$-complete \cite{MR1567289}.

 {\em Monotone Not-All-Equal 3-Sat.}\\
\textsc{Instance}:
Set $X$ of variables, collection $C$ of clauses over $X$ such that each
clause $c \in C$ has $\mid c  \mid = 3$ and there is no negation in the formula.\\
\textsc{Question}: Is there a truth assignment for $X$ such that each clause in $C$ has at
least one true literal and at least one false literal?\\

We show that this problem remains $ \mathbf{NP}$-complete when every variable has exactly three occurrences and every clause contains two or three variables.

 {\em Cubic Monotone NAE (2,3)-Sat.}\\
\textsc{Instance}: Set $X$ of variables, collection $C$ of clauses over $X$ such that each
clause $c \in C$ has $\mid c  \mid \in \{ 2, 3\}$, every variable appears in
exactly three clauses and there is no negation in the formula.\\
\textsc{Question}: Is there a truth assignment for $X$ such that each clause in $C$ has at
least one true literal and at least one false literal?\\

\begin{figure}[ht]
\begin{center}
\includegraphics[scale=.6]{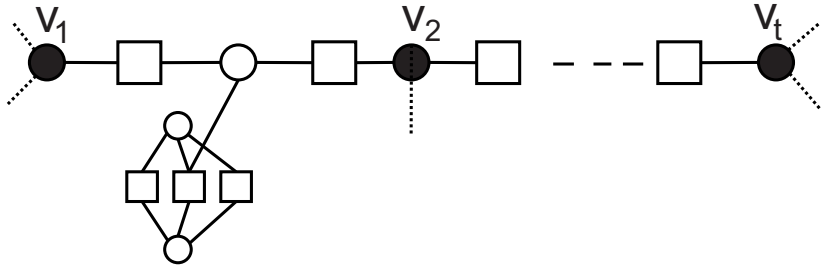}
\caption{The chain  $Ch_t$  forces black variables to be equal.
} \label{graphB}
\end{center}
\end{figure}

Consider the graph $Ch_t$ which is shown in Figure \ref{graphB}. Every circle is a variable, every square is a clause and the clause $c$ contains the variable $x$ if and only if they are connected in the graph. Also, every variable occurs at most three times and each clause contains two or three variables. The chain  $Ch_t$  forces black variables to be equal.

Let $\Psi$ be an instance of {\em Monotone NAE 3-Sat}.
For each variable $v$ that occurs
in $k > 3$ clauses, we replace it with the chain $Ch_{k-2}$.
Also, for every variable $v$ that occurs
in $2 $ clauses, create two new variables $v_1$ and $v_2$, next add three clauses $(v_1 \vee v_2)$, $(v_1 \vee v_2)$ and
$(v_1 \vee v_2 \vee v)$. For every variable $v$ that occurs
in one clause, perform a similar method.
Therefore, we obtain a formula $ \Phi$ which is an instance of
 {\em Cubic Monotone NAE (2,3)-Sat} and also $ \Phi$
has an NAE truth assignment if and
only if $\Psi$ has an NAE truth assignment.

{{\bf Step 2.}} Now we reduce {\em Cubic Monotone NAE (2,3)-Sat } to
our problem in polynomial time.
Consider an instance $ \Phi $, we transform this into a  graph $G_{\Phi}$
such that $ reg(G_{\Phi})=2$ if and only if $\Phi$ has an  NAE truth assignment.
 We use two gadgets $ H_c$ and $I_c$ which are shown in
Figure \ref{graphA}. The graph $G_{\Phi}$ has
a copy of $H_c $ for each clause $c \in C$ with $\mid c \mid =3$, a copy of $I_c $ for each clause $c \in C$ with $\mid c \mid =2$ and an  vertex $x$ for each variable $x \in X$.
For each clause $c =y \vee z \vee w$, where  $y,w,z \in X  $ add the edges $u_c y $, $u_c z  $ and $u_c w $. Also, for each clause $c =y \vee z $, where  $y,w \in X  $ add the edges $v_c y $ and $v_c z  $. The degree of every vertex in $G$ is $3$ or $6$ and $G_{\Phi}$ is $3$-colorable.

\begin{figure}[ht]
\begin{center}
\includegraphics[scale=.6]{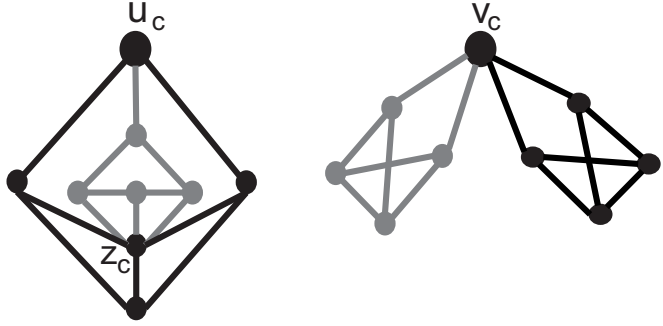}
\caption{The two gadgets $ H_c$ and $I_c$. $I_c$ is on the right hand side of the figure.
} \label{graphA}
\end{center}
\end{figure}

First, assume that $reg(G_{\Phi})=2$ and   $G_1,  G_2$ is a regular partition of $G$ such that $G_i$ is $(r_i)$-regular.
Since $G_{\Phi}$ has vertices with degrees $3$ and $6$, so $r_1=r_2=3$.
For every $x \in X$, since $x$ has degree 3, so all edges incident with $x$ are in one part. For every $x\in X$, if all edges incident with $x$ are in $G_1$, put $\Gamma(x)=true$ and if all edges incident with $x$ are in $G_2$, put $\Gamma(x)=false$.
According to the construction of $H_c$ and $I_c$, in each of them the black edges are in one part and the gray edges are in another part.
Therefore, for every clause $c=y \vee z \vee w$, at most two of the three edges $u_c y $, $u_c z  $ and $u_c w $ are in $G_1$. Also, at most two of the three edges $u_c y $, $u_c z  $ and $u_c w $ are in $G_2$.
By a similar argument, for every clause $c=y \vee z $, exactly one of the two edges $v_c y $ and $v_c z  $ is in $G_1$.
Hence, $\Gamma$ is an NAE satisfying assignment.
On the other hand, suppose that $\Phi$ is NAE satisfiable with the satisfying assignment $\Gamma : X  \rightarrow \{true, false\}$. For every variable $x\in X$, put all edges incident with $x$ in $G_1$ if and only if $\Gamma (x)=true$, also for every subgraph $I_c$, put the black edges in $G_1$ and the gray edges in $G_2$. Finally, for every clause $c=y \vee z \vee w$, if $c$ has two true literals, in the subgraph $H_c$, put the gray edges in $G_1$ and the black edges in $G_2$ and otherwise put the black edges in $G_1$ and the gray edges in $G_2$.
This completes the proof.

}\end{alii}

A. Ganesan et al.  in \cite{kulli2} asked to determine whether  $reg(G)\leq \Delta(G)$ holds for all
connected graphs $G$.
We show that not only  there exists a counterexample for the above bound but also for a given connected graph $G$ deciding whether $reg(G)=\Delta(G)+1$ is
{\bf NP}-complete.

\begin{aliii}{
For a cubic graph $G$, define $\sigma(G)$ to be the minimum
number such that there exists a $4$-edge-coloring of $G$ with
$\sigma(G)$ edges assigned the fourth color.
Kochol et al.  proved that for a given  cyclically 6-edge-connected 3-regular graph $G$
approximating $\sigma(G)$ with an error $O(n^{1-\varepsilon})$ is $ \mathbf{NP}$-hard, where $n$ denotes the number of vertices of the graph $G$ and $0 < \varepsilon < 1$ is a constant  \cite{kochol}.
In other words, they proved that the problem to
decide whether $ \sigma(G)\in [0,n^{1-\varepsilon}]$ is $ \mathbf{NP}$-complete in the class of cyclically
6-edge-connected cubic graphs
(note that a graph is cyclically $k$-edge-connected
if deleting fewer than $k$ edges does not result in a graph having at least two
components containing cycles).

In fact, kochol et al. proved that, the problem to decide whether a cyclically 6-edge-connected
cubic graph has a 3-edge-coloring is $ \mathbf{NP}$-complete and
for a given cyclically $6$-edge-connected cubic graph $H$ with $m$ vertices, in polynomial time they construct a graph $G$  which is
cyclically $6$-edge connected cubic graph with $ n$ vertices such that:\\ \\
(\textsc{i}): $H$ has a $3$-edge-coloring if and only if $G$ has $3$-edge-coloring.\\
(\textsc{ii}): If $G$ is not 3-edge-colorable, then $\sigma(G)>  n^{1-\varepsilon} $.

Therefore, Kochol et al. construct a family of cubic graphs such that every graph $G$ in this family has a 3-edge-coloring or
 $\sigma(G)>  n^{1-\varepsilon}$, where $n$ is the number of vertices of $G$  \cite{kochol}. Also, they show that determining the 3-edge-colorability in this family is $ \mathbf{NP}$-complete  \cite{kochol}. Consequently, by this reduction the following promise problem is  $ \mathbf{NP}$-complete.

{\em  Problem $\Theta$.}\\
\textsc{Instance}: A cyclically 6-edge-connected cubic graph $G$.\\
\textsc{Question}: Is $\chi'(G)\leq 3$?\\
\textsc{Promise}: $\sigma(G)\neq 1$.\\

We reduce Problem $\Theta$ to our problem, let $G$ be  a given cyclically 6-edge-connected 3-regular graph and $e=uv \in E(G)$ be an arbitrary edge, remove $e$ from $E(G)$ and add three vertices $v'$, $v''$ and $v'''$ to $V(G)$; join $v'$ to $v$,$v''$ and $v'''$. Call the resulting graph $G'$. Since $\sigma(G)\neq 1$, therefore if $G$ is not 3-edge-colorable
then, $G\setminus e$ is also not 3-edge-colorable, so if $G$ is not 3-edge-colorable, then $G'$ is not 3-edge-colorable.
Consequently, since $d(v'')=d(v''')=1$ and $d(v')=3$ we have $reg(G')=3$ if and only if $\sigma(G) =0 $. This completes the proof.
}\end{aliii}

\section{Acknowledgment}
\label{}

The authors would like to thank the anonymous referee for his/her useful comments and suggestions, which helped to
improve the presentation of this paper.

\bibliographystyle{plain}
\bibliography{luckyref}

\begin{thebibliography}{10}

\bibitem{MR2510539}
Krzysztof Bry{\'s} and Zbigniew Lonc.
\newblock Polynomial cases of graph decomposition: a complete solution of
  {H}olyer's problem.
\newblock {\em Discrete Math.}, 309(6):1294--1326, 2009.

\bibitem{MR1460720}
Dorit Dor and Michael Tarsi.
\newblock Graph decomposition is {NP}-complete: a complete proof of {H}olyer's
  conjecture.
\newblock {\em SIAM J. Comput.}, 26(4):1166--1187, 1997.

\bibitem{kulli2}
Ashwin Ganesan and Radha~R. Iyer.
\newblock The regular number of a graph.
\newblock {\em Journal of Discrete Mathematical Sciences and Cryptography},
  15(2 -3):149--157, 2012.

\bibitem{MR1567289}
M.~R. Garey and D.~S. Johnson.
\newblock {\em {C}omputers and intractability: {A} guide to the theory of
  {$NP$}-completeness}.
\newblock W. H. Freeman, San Francisco, 1979.

\bibitem{MR635429}
Ian Holyer.
\newblock The {NP}-completeness of some edge-partition problems.
\newblock {\em SIAM J. Comput.}, 10(4):713--717, 1981.

\bibitem{kochol}
Martin Kochol, Nad'a Krivo{\v{n}}{\'a}kov{\'a}, Silvia Smejov{\'a}, and
  Katar{\'{\i}}na {\v{S}}rankov{\'a}.
\newblock Complexity of approximation of 3-edge-coloring of graphs.
\newblock {\em Inform. Process. Lett.}, 108(4):238--241, 2008.

\bibitem{kulli}
V.~R. Kulli, B.~Janakiram, and R.~R. Iyer.
\newblock Regular number of a graph.
\newblock {\em Journal of Discrete Mathematical Sciences and Cryptography},
  4(1):57--64, 2001.

\bibitem{MR689264}
Daniel Leven and Zvi Galil.
\newblock N{P}-completeness of finding the chromatic index of regular graphs.
\newblock {\em J. Algorithms}, 4(1):35--44, 1983.

\bibitem{MR1403227}
Zbigniew Lonc.
\newblock On the complexity of some edge-partition problems for graphs.
\newblock {\em Discrete Appl. Math.}, 70(2):177--183, 1996.

\bibitem{MR679633}
B.~P{\'e}roche.
\newblock Complexity of the linear arboricity of a graph.
\newblock {\em RAIRO Rech. Op\'er.}, 16(2):125--129, 1982(in French).

\bibitem{MR2297163}
Michael~D. Plummer.
\newblock Graph factors and factorization: 1985--2003: a survey.
\newblock {\em Discrete Math.}, 307(7-8):791--821, 2007.

\bibitem{MR0180505}
V.~G. Vizing.
\newblock On an estimate of the chromatic class of a {$p$}-graph.
\newblock {\em Diskret. Analiz No.}, 3:25--30, 1964.

\bibitem{MR1367739}
Douglas~B. West.
\newblock {\em Introduction to graph theory}.
\newblock Prentice Hall Inc., Upper Saddle River, NJ, 1996.

\end{thebibliography}

\end{document}